\documentclass[a4paper,11pt]{amsart}
\usepackage{graphicx}
 \usepackage{mathptmx}
\usepackage{amsmath}
\usepackage{amssymb}
\usepackage{enumitem}
\usepackage{xcolor}

\newmuskip\pFqmuskip

\newcommand*\pFq[6][8]{%
  \begingroup 
  \pFqmuskip=#1mu\relax
  \mathcode`=\string"8000
  \begingroup\lccode`\~=`\,
  \lowercase{\endgroup\let~}\pFqcomma
  F^{#2}_{#3}{\left(\genfrac..{0pt}{}{#4}{#5}\bigg|#6\right)}%
  \endgroup
}
\newcommand{\pFqcomma}{\mskip\pFqmuskip}

\newtheorem{theorem}{Theorem}
\newtheorem{lemma}[theorem]{Lemma}
\newtheorem{corollary}[theorem]{Corollary}

\begin{document}

\title[A note on new type degenerate Bernoulli numbers]{A note on new type degenerate Bernoulli numbers}

\author{Taekyun  Kim}
\address{Department of Mathematics, Kwangwoon University, Seoul 139-701, Republic of Korea}
\email{tkkim@kw.ac.kr}

\author{Dae San  Kim }
\address{Department of Mathematics, Sogang University, Seoul 121-742, Republic of Korea}
\email{dskim@sogang.ac.kr}


\maketitle

\begin{abstract}
Studying degenerate versions of various special polynomials have become an active area of research and yielded many interesting arithmetic and combinatorial results.
Here we introduce a degenerate version of polylogarithm function, called the degenerate polylogarithm function. Then we construct new type degenerate Bernoulli polynomials and numbers, called degenerate poly-Bernoulli polynomials and numbers, by using the degenerate polylogarithm function and derive several properties on the degenerate poly-Bernoulli numbers.
\end{abstract}

\section{Introduction}
As is well known, for $s\in\mathbb{C}$, the polylogarithm function is defined by a power series in $z$, which is also a Dirichlet series in $s$:
\begin{equation}
	\mathrm{Li}_{s}(z)=\sum_{n=1}^{\infty}\frac{z^{n}}{n^{s}}=z+\frac{z^{2}}{2^{s}}+\frac{z^{3}}{3^{s}}+\cdots,\quad(\mathrm{see}\ [6,9,22]). \label{1}
\end{equation}
This definition is valid for arbitrary complex order $s$ and for all complex arguments $z$ with $|z|<1$: it can be extended to $|z|\ge 1$ by analytic continuation. \\
From \eqref{1}, we note that
\begin{equation}
	\mathrm{Li}_{1}(z)=\sum_{n=1}^{\infty}\frac{z^{n}}{n}=-\log(1-z). \label{2}
\end{equation}
For any nonzero $\lambda\in\mathbb{R}$ (or $\mathbb{C}$), the degenerate exponential function is defined by
\begin{equation}
	e_{\lambda}^{x}(t)=(1+\lambda t)^{\frac{x}{\lambda}},\quad e_{\lambda}(t)=(1+\lambda t)^{\frac{1}{\lambda}}=e_{\lambda}^{1}(t),\quad (\mathrm{see}\ [1,14,15,17]). \label{3}
\end{equation}
By Taylor expansion, we get
\begin{equation}
	e_{\lambda}^{x}(t)=\sum_{n=0}^{\infty}(x)_{n,\lambda}\frac{t^{n}}{n!},\quad(\mathrm{see}\ [13,14,15,16,17,18]),\label{4}
\end{equation}
where $\displaystyle (x)_{0,\lambda}=1,\ (x)_{n,\lambda}=(x-\lambda)(x-2\lambda)\cdots(x-(n-1)\lambda),\ (n\ge 1)\displaystyle$. \\
Note that
\begin{displaymath}
	\lim_{\lambda\rightarrow 0}e_{\lambda}^{x}(t)=\sum_{n=0}^{\infty}\frac{x^{n}t^{n}}{n!}=e^{xt}.
\end{displaymath}
In [1,2], Carlitz considered the degenerate Bernoulli polynomials given by
\begin{equation}
	\frac{t}{e_{\lambda}(t)-1}e_{\lambda}^{x}(t)=\frac{t}{(1+\lambda t)^{\frac{1}{\lambda}}-1}(1+\lambda t)^{\frac{x}{\lambda}}=\sum_{n=0}^{\infty}\beta_{n,\lambda}(x)\frac{t^{n}}{n!}.\label{5}
\end{equation}
When $x=0$, $\beta_{n,\lambda}=\beta_{n,\lambda}(0)$ are called the degenerate Bernoulli numbers. Note that
\begin{displaymath}
	\lim_{\lambda\rightarrow 0}\beta_{n,\lambda}(x)=B_{n}(x),\quad (n\ge 0),
\end{displaymath}
where $B_{n}(x)$ are the ordinary Bernoulli polynomials given by
\begin{equation}
	\frac{t}{e^{t}-1}e^{xt}=\sum_{n=0}^{\infty}B_{n}(x)\frac{t^{n}}{n!},\quad(\mathrm{see}\ [1-26]).\label{6}
\end{equation}
The Stirling numbers of the second kind are defined as
\begin{equation}
	x^{n}=\sum_{l=0}^{n}S_{2}(n,l)(x)_{l},\quad(n\ge 0),\quad (\mathrm{see}\ [10,12,14,20,22]),\label{7}
\end{equation}
where $(x)_{0}=1,\ (x)_{n}=x(x-1)(x-2)\cdots(x-n+1),\ (n\ge 1)$. \\
Thus, we easily get
\begin{equation}
	\frac{1}{k!}(e^{t}-1)^{k}=\sum_{n=k}^{\infty}S_{2}(n,k)\frac{t^{n}}{n!},\quad (k\ge 0). \label{8}
\end{equation}
In [11], Kim considered the degenerate Stirling numbers of the second kind which are defined as
\begin{equation}
	(x)_{n,\lambda}=\sum_{k=0}^{n}S_{2,\lambda}(n,k)(x)_{k},\quad (n\ge 0)\label{9}.
\end{equation}
Note that
\begin{displaymath}
	\lim_{\lambda\rightarrow 0}S_{2,\lambda}(n,k)=S_{2}(n,k).
\end{displaymath}
The generating function of the degenerate Stirling numbers of the second kind is given by
\begin{equation}
	\frac{1}{k!}(e_{\lambda}(t)-1)^{k}=\sum_{n=k}^{\infty}S_{2,\lambda}(n,k)\frac{t^{n}}{n!},\quad(k\ge 0),\quad(\mathrm{see}\ [11]). \label{10}
\end{equation}
In this paper, we will introduce the degenerate polylogarithm function as a degenerate version of the polylogarithm function for $s=k \in \mathbb{Z}$.
Then we will construct new type degenerate Bernoulli polynomials and numbers, called degenerate poly-Bernoulli polynomials and numbers, by using the degenerate polylogarithm function and derive several properties on the degenerate poly-Bernoulli numbers.

\section{New type degenerate Bernoulli numbers and polynomials}
We define the degenerate logarithm function $\log_{\lambda}(1+t)$, which is the inverse of the degenerate exponential function $e_{\lambda}(t)$ and the motivation for the definition of degenerate polylogarithm function, as:
\begin{equation}
	\log_{\lambda}(1+t)=\sum_{n=1}^{\infty}\lambda^{n-1}(1)_{n,1/\lambda}\frac{t^{n}}{n!}.\label{11}
\end{equation}
From \eqref{4} and \eqref{11}, we note that
\begin{align}
	\log_{\lambda}(1+t)&=\sum_{n=1}^{\infty}\lambda^{n-1}\bigg(1-\frac{1}{\lambda}\bigg)\bigg(1-\frac{2}{\lambda}\bigg)\cdots\bigg(1-\frac{1}{\lambda}(n-1)\bigg)\frac{t^{n}}{n!} \nonumber \label{12}\\
	&=\frac{1}{\lambda}\sum_{n=1}^{\infty}(\lambda)_n\frac{t^{n}}{n!} \nonumber \\
	&=\frac{1}{\lambda}\big((1+t)^{\lambda}-1\big). \nonumber
\end{align}
\begin{lemma}
	For $\lambda\in\mathbb{R}$, we have
	\begin{displaymath}
		\log_{\lambda}(1+t)=\frac{1}{\lambda}\big((1+t)^{\lambda}-1\big).
	\end{displaymath}
In addition, $e_{\lambda}\big(\log_{\lambda}(1+t)\big)=1+t$.
\end{lemma}
It is easy to show that
\begin{displaymath}
	\lim_{\lambda\rightarrow 0}\log_{\lambda}(1+t)=\sum_{n=1}^{\infty}(-1)^{n-1}\frac{t^{n}}{n}=\log(1+t).
\end{displaymath}
For $k\in\mathbb{Z}$, we define the degenerate polylogarithm function as
\begin{equation}
	l_{k,\lambda}(x)=\sum_{n=1}^{\infty}\frac{(-\lambda)^{n-1}(1)_{n,1/\lambda}}{(n-1)!n^{k}}x^{n},\quad (|x|<1). \label{13}
\end{equation}
Note that
\begin{displaymath}
	\lim_{\lambda\rightarrow 0}l_{k,\lambda}(x)=\sum_{n=1}^{\infty}\frac{x^{n}}{n^{k}}=\mathrm{Li}_{k}(x),\quad (\mathrm{see}\ [3,22]).
\end{displaymath}
From \eqref{13}, we note that
\begin{equation}
	\frac{d}{dx}l_{k,\lambda}(x)=\frac{d}{dx}\sum_{n=1}^{\infty}\frac{(-\lambda)^{n-1}(1)_{n,1/\lambda}}{(n-1)!n^{k}}x^{n}=\frac{1}{x}l_{k-1,\lambda}(x)\label{14} .
\end{equation}
For $k\ge 2$, \eqref{14} can be written in the form of iterated integral which is given by
\begin{equation}
	l_{k,\lambda}(x)=\int_{0}^{x}\underbrace{\frac{1}{t}\int_{0}^{t}\frac{1}{t}\int_{0}^{t}\frac{1}{t}\cdots\int_{0}^{t}}_{(k-2)-\mathrm{times}}\frac{1}{t}l_{1,\lambda}(t)dtdt\cdots dt \label{15}
\end{equation}
By \eqref{11} and \eqref{13}, we get
\begin{equation}
	l_{1,\lambda}(x)=\sum_{n=1}^{\infty}\frac{(-\lambda)^{n-1}(1)_{n,1/\lambda}}{n!}x^{n}=-\log_{\lambda}(1-x)\label{16} .
\end{equation}
Thus, by \eqref{15} and \eqref{16}, for $k \ge 2$ we get
\begin{equation}
	l_{k,\lambda}(x)=-\int_{0}^{x}\underbrace{\frac{1}{t}\int_{0}^{t}\frac{1}{t}\int_{0}^{t}\frac{1}{t}\cdots\int_{0}^{t}}_{(k-2)-\mathrm{times}}\frac{1}{t}\log_{\lambda}(1-t)dtdt\cdots dt. \label{17}
\end{equation}
For $k\in\mathbb{Z}$, we define the new type degenerate Bernoulli numbers, which are called the degenerate poly-Bernoulli numbers, as
\begin{equation}
	\frac{1}{x}l_{k,\lambda}(x)\bigg|_{x=1-e_{\lambda}(-t)}=\frac{1}{1-e_{\lambda}(-t)}l_{k,\lambda}(1-e_{\lambda}(-t))=\sum_{n=0}^{\infty}\beta_{n,\lambda}^{(k)}\frac{t^{n}}{n!}. \label{18}
\end{equation}
 Note that
 \begin{equation}
 	\sum_{n=0}^{\infty}\beta_{n,\lambda}^{(1)}\frac{t^{n}}{n!}=\frac{1}{1-e_{\lambda}(-t)}l_{1,\lambda}(1-e_{\lambda}(-t))=\frac{-t}{e_{\lambda}(-t)-1}=\sum_{n=0}^{\infty}\frac{(-1)^{n}\beta_{n,\lambda}}{n!}t^{n}. \label{19}
 \end{equation}
 Comparing the coefficients on both sides of \eqref{19}, we have
 \begin{displaymath}
 	\beta_{n,\lambda}^{(1)}=(-1)^{n}\beta_{n,\lambda},\quad (n\ge 0).
 \end{displaymath}
 Now, we consider the new type degenerate Bernoulli polynomials which are called the degenerate poly-Bernoulli polynomials and given by
\begin{equation}
	\frac{l_{k,\lambda}(1-e_{\lambda}(-t))}{1-e_{\lambda}(-t)}e_{\lambda}^{x}(-t)=\sum_{n=0}^{\infty}\beta_{n,\lambda}^{(k)}(x)\frac{t^{n}}{n!}. \label{20}
\end{equation}
Note here that $\beta_{n,\lambda}^{(k)}=\beta_{n,\lambda}^{(k)}(0)$.
From \eqref{18} and \eqref{20}, we note that
\begin{align}
	\sum_{n=0}^{\infty}\beta_{n,\lambda}^{(k)}(x)\frac{t^{n}}{n!}&=\frac{l_{k,\lambda}(1-e_{\lambda}(-t))}{1-e_{\lambda}(-t)}e_{\lambda}^{x}(-t)\label{21}\\
	&=\sum_{l=0}^{\infty}\beta_{l,\lambda}^{(k)}\frac{t^{l}}{l!}\sum_{m=0}^{\infty}\frac{(-1)^{m}(x)_{m,\lambda}}{m!}t^{m}\nonumber\\
	&=\sum_{n=0}^{\infty}\bigg(\sum_{l=0}^{n}\binom{n}{l}\beta_{l,\lambda}^{(k)}(-1)^{n-l}(x)_{n-l,\lambda}\bigg)\frac{t^{n}}{n!}\nonumber
\end{align}
Comparing the coefficients on both sides of \eqref{21}, we have
\begin{equation}
	\beta_{n,\lambda}^{(k)}(x)=\sum_{l=0}^{n}\binom{n}{l}\beta_{l,\lambda}^{(k)}(-1)^{n-l}(x)_{n-l,\lambda},\quad (n\ge 0). \label{22}
\end{equation}
Now, we observe that
\begin{align}
	\frac{d}{dx}e_{\lambda}(-x)&=\frac{d}{dx}\sum_{l=0}^{\infty}\frac{(-1)^{l}(1)_{l,\lambda}}{l!}x^{l}=\sum_{l=1}^{\infty}\frac{(-1)^{l}(1)_{l,\lambda}}{(l-1)!}x^{l-1}\label{23}\\
	&=-\sum_{l=0}^{\infty}\frac{(-1)^{l}(1)_{l+1,\lambda}}{l!}x^{l}=-\sum_{l=0}^{\infty}\frac{(-1)^{l}(1)_{l,\lambda}}{l!}x^{l}(1-l\lambda) \nonumber \\
	&=-e_{\lambda}(-x)+\lambda\sum_{l=1}^{\infty}\frac{(-1)^{l}(1)_{l,\lambda}}{(l-1)!}x^{l}=-e_{\lambda}(-x)+\lambda x\frac{d}{dx}e_{\lambda}(-x). \nonumber
\end{align}
Thus, by \eqref{23}, we get
\begin{equation}
	(1-\lambda x)\frac{d}{dx}e_{\lambda}(-x)=-e_{\lambda}(-x).\label{24}
\end{equation}
Therefore, by \eqref{24}, we obtain the following lemma.
\begin{lemma}
	For $\lambda\in\mathbb{R}$, we have
	\begin{displaymath}
		\frac{d}{dx}e_{\lambda}(-x)=-\frac{1}{1-\lambda x}e_{\lambda}(-x)=-e_{\lambda}^{1-\lambda}(-x).
	\end{displaymath}
\end{lemma}
By Lemma 2, we easily get
\begin{equation}
	\frac{d}{dx}(1-e_{\lambda}(-x))=\frac{1}{1-\lambda x}e_{\lambda}(-x)=e_{\lambda}^{1-\lambda}(-x)\label{25}
\end{equation}
From \eqref{14}, \eqref{16}, \eqref{18} and \eqref{25}, for $k \ge 2$ we have
\begin{align}
&\sum_{n=0}^{\infty}\beta_{n,\lambda}^{(k)}\frac{x^{n}}{n!}=\frac{1}{1-e_{\lambda}(-x)}l_{k,\lambda}(1-e_{\lambda}(-x))\label{26}\\
&=\frac{1}{1-e_{\lambda}(-x)}\int_{0}^{x}\underbrace{\frac{e_{\lambda}^{1-\lambda}(-t)}{1-e_{\lambda}(-t)}\int_{0}^{t}\frac{e_{\lambda}^{1-\lambda}(-t)}{1-e_{\lambda}(-t)}\cdots\int_{0}^{t}}_{(k-2)-\mathrm{times}}\frac{e_{\lambda}^{1-\lambda}(-t)}{1-e_{\lambda}(-t)}t \,\,dtdt\cdots dt.\nonumber
\end{align}
Therefore, by \eqref{26}, we obtain the following theorem.
\begin{theorem}
For $k \ge 2$, we have
\begin{displaymath}
\sum_{n=0}^{\infty}\beta_{n,\lambda}^{(k)}\frac{x^{n}}{n!}=\frac{1}{1-e_{\lambda}(-x)}\int_{0}^{x}\underbrace{\frac{e_{\lambda}^{1-\lambda}(-t)}{1-e_{\lambda}(-t)}\int_{0}^{t}\frac{e_{\lambda}^{1-\lambda}(-t)}{1-e_{\lambda}(-t)}\cdots\int_{0}^{t}}_{(k-2)-\mathrm{times}}\frac{e_{\lambda}^{1-\lambda}(-t)}{1-e_{\lambda}(-t)} t\,\,dtdt\cdots dt.
\end{displaymath}
\end{theorem}
From \eqref{18}, we can derive the following equation:
\begin{align}
	\sum_{n=0}^{\infty}\beta_{n,\lambda}^{(k)}\frac{t^{n}}{n!}&=\sum_{n=1}^{\infty}\frac{(-\lambda)^{n-1}(1)_{n,1/\lambda}}{(n-1)!n^{k}}\big(1-e_{\lambda}(-t)\big)^{n-1}\label{27} \\
	&=\sum_{m=0}^{\infty}\frac{(-\lambda)^{m}(1)_{m+1,1/\lambda}}{m!(m+1)^{k}}(1-e_{\lambda}(-t))^{m}\nonumber\\
	&=\sum_{m=0}^{\infty}\frac{(-\lambda)^{m}}{(m+1)^{k}}(1)_{m+1,1/\lambda}\sum_{n=m}^{\infty}(-1)^{m-n}S_{2,\lambda}(n,m)\frac{t^{n}}{n!} \nonumber\\
	&=\sum_{n=0}^{\infty}\bigg((-1)^{n}\sum_{m=0}^{n}\frac{\lambda^{m}(1)_{m+1,1/\lambda}}{(m+1)^{k}}S_{2,\lambda}(n,m)\bigg)\frac{t^{n}}{n!}. \nonumber
\end{align}
Therefore, by comparing the coefficients on both sides of \eqref{27}, we obtain the following theorem.
\begin{theorem}
	For $n\ge 0$, we have
	\begin{displaymath}
		\beta_{n,\lambda}^{(k)}=(-1)^{n}\sum_{m=0}^{n}\frac{\lambda^{m}(1)_{m+1,1/\lambda}}{(m+1)^{k}}S_{2,\lambda}(n,m).
	\end{displaymath}
\end{theorem}
Note that
\begin{displaymath}
	(-1)^{n}B_{n}=\lim_{\lambda\rightarrow 0}\beta_{n,\lambda}^{(1)}=(-1)^n\sum_{m=0}^{n}\frac{m!}{m+1}(-1)^{m}S_{2}(n,m),\quad (n\ge 0).
\end{displaymath}
For $k=2$, by Theorem 3, we get
\begin{align}
	\sum_{n=0}^{\infty}\beta_{n,\lambda}^{(2)}\frac{x^{n}}{n!}&=\frac{1}{1-e_{\lambda}(-x)}\int_{0}^{x}\frac{t}{1-e_{\lambda}(-t)}e_{\lambda}^{1-\lambda}(-t)dt \label{28} \\
	&=\frac{1}{1-e_{\lambda}(-x)}\int_{0}^{x}\sum_{n=0}^{\infty}\beta_{n,\lambda}(1-\lambda)(-1)^{n}\frac{t^{n}}{n!}dt \nonumber \\
	&=\frac{x}{1-e_{\lambda}(-x)}\sum_{l=0}^{\infty}\frac{\beta_{l,\lambda}(1-\lambda)}{l+1}\frac{(-1)^lx^{l}}{l!}\nonumber \\
	&=\sum_{m=0}^{\infty}(-1)^{m}\beta_{m,\lambda}\frac{x^{m}}{m!}\sum_{l=0}^{\infty}(-1)^l\frac{\beta_{l,\lambda}(1-\lambda)}{l+1}\frac{x^{l}}{l!}\nonumber \\
	&=\sum_{n=0}^{\infty}\bigg((-1)^{n}\sum_{m=0}^{n}\binom{n}{m}\beta_{m,\lambda}\frac{\beta_{n-m,\lambda}(1-\lambda)}{n-m+1}\bigg)\frac{x^{n}}{n!}. \nonumber
\end{align}
Therefore, by comparing the coefficients on both sides of \eqref{28}, we obtain the following theorem.
\begin{theorem}
	For $n\ge 0$, we have
	\begin{displaymath}
		\beta_{n,\lambda}^{(2)}=(-1)^{n}\sum_{m=0}^{n}\binom{n}{m}\beta_{m,\lambda}\frac{\beta_{n-m,\lambda}(1-\lambda)}{n-m+1}=(-1)^{n}\sum_{m=0}^{n}\binom{n}{m}\beta_{n-m,\lambda}\frac{\beta_{m,\lambda}(1-\lambda)}{m+1}.
	\end{displaymath}
\end{theorem}
In general, from (25), we note that
\begin{equation}\label{28-1}
\sum_{n=0}^{\infty}\beta_{n,\lambda}^{(k)}\frac{x^{n}}{n!}=\frac{1}{1-e_{\lambda}(-x)}\int_{0}^{x}\frac{e_{\lambda}^{1-\lambda}(-t)}{1-e_{\lambda}(-t)}\int_{0}^{t}\frac{e_{\lambda}^{1-\lambda}(-t)}{1-e_{\lambda}(-t)}\cdots\int_{0}^{t}\frac{e_{\lambda}^{1-\lambda}(-t)}{1-e_{\lambda}(-t)} t\,\,dtdt\cdots dt
\end{equation}
\begin{align*}
	&=\sum_{n_{1},n_{2},\cdots,n_{k-1}=0}^{\infty}\frac{1}{n_{1}!n_{2}!\cdots n_{k-1}!} \frac{\beta_{n_{1},\lambda}(1-\lambda)}{n_{1}+1}\frac{\beta_{n_{2},\lambda}(1-\lambda)}{n_{1}+n_{2}+1}\cdots\frac{\beta_{n_{k-1},\lambda}(1-\lambda)}{n_{1}+\cdots+n_{k-1}+1}(-x)^{n_{1}+\cdots+n_{k-1}}\frac{x}{1-e_{\lambda}(-x)}\\
	&=\quad\quad\sum_{n=0}^{\infty}(-1)^{n}\sum_{n_{1}+\cdots+n_{k}=n}\binom{n}{n_{1},n_{2},\dots,n_{k}}\frac{\beta_{n_{1},\lambda}(1-\lambda)}{n_{1}+1}\frac{\beta_{n_{2},\lambda}(1-\lambda)}{n_{1}+n_{2}+1}\cdots\frac{\beta_{n_{k-1},\lambda}(1-\lambda)}{n_{1}+n_{2}+\cdots+n_{k-1}+1}\beta_{n_{k},\lambda}\frac{x^{n}}{n!}.
\end{align*}
Therefore, by comparing the coefficients on both sides of \eqref{28-1}, we obtain the following theorem.
\begin{theorem}
	For $k\in\mathbb{N}$ and $n\ge 0$, we have
\begin{displaymath}
	\beta_{n,\lambda}^{(k)}=(-1)^{n}\sum_{n_{1}+\cdots+n_{k}=n}\binom{n}{n_{1},n_{2},\dots,n_{k}}\frac{\beta_{n_{1},\lambda}(1-\lambda)}{n_{1}+1}\frac{\beta_{n_{2},\lambda}(1-\lambda)}{n_{1}+n_{2}+1}\cdots\frac{\beta_{n_{k-1},\lambda}(1-\lambda)}{n_{1}+n_{2}+\cdots+n_{k-1}+1}\beta_{n_{k},\lambda}.
\end{displaymath}
\end{theorem}

From \eqref{18}, we observe that
\begin{align}
	l_{k,\lambda}\big(1-e_{\lambda}(-t)\big)&=\big(1-e_{\lambda}(-t)\big) \sum_{l=0}^{\infty}\beta_{l,\lambda}^{(k)}\frac{t^{l}}{l!}\label{34} \\
	&=\bigg(1-\sum_{m=0}^{\infty}\frac{(-1)^{m}(1)_{m,\lambda}}{m!}t^{m}\bigg)\sum_{l=0}^{\infty}\beta_{l,\lambda}^{(k)}\frac{t^{l}}{l!}\nonumber \\
	&=\sum_{n=0}^{\infty}\bigg(\beta_{n,\lambda}^{(k)}-\sum_{l=0}^{n}\binom{n}{l}\beta_{l,\lambda}^{(k)}(-1)^{n-l}(1)_{n-l,\lambda}\bigg)\frac{t^{n}}{n!} \nonumber \\
&=\sum_{n=1}^{\infty}\big(\beta_{n,\lambda}^{(k)}-\beta_{n,\lambda}^{(k)}(1)\big)\frac{t^{n}}{n!}. \nonumber
\end{align}
On the other hand,
\begin{align}
	l_{k,\lambda}\big(1-e_{\lambda}(-t)\big)&=\sum_{m=1}^{\infty}\frac{(1)_{m,1/\lambda}(-\lambda)^{m-1}}{(m-1)!m^{k}}\big(1-e_{\lambda}(-t)\big)^{m} \label{35} \\
	&=\sum_{m=1}^{\infty}\frac{(1)_{m,1/\lambda}}{m^{k-1}}(-\lambda)^{m-1}\frac{1}{m!}\big(1-e_{\lambda}(-t)\big)^{m}\nonumber\\
	&=\sum_{m=1}^{\infty}\frac{(1)_{m,1/\lambda}(-\lambda)^{m-1}}{m^{k-1}}\sum_{n=m}^{\infty}S_{2,\lambda}(n,m)(-1)^{n-m}\frac{t^{n}}{n!} \nonumber \\
	&=\sum_{n=1}^{\infty}\bigg(\sum_{m=1}^{n}\frac{(1)_{m,1/\lambda}(-1)^{n-1}}{m^{k-1}}\lambda^{m-1}S_{2,\lambda}(n,m)\bigg)\frac{t^{n}}{n!}. \nonumber
\end{align}
Therefore, by \eqref{34} and \eqref{35}, we obtain the following theorem.
\begin{theorem}
	For $k\in\mathbb{Z}$, we have
	\begin{displaymath}
		\quad \beta_{n,\lambda}^{(k)}(1)-\beta_{n,\lambda}^{(k)}=(-1)^{n}\sum_{m=1}^{n}\frac{(1)_{m,1/\lambda}}{m^{k-1}}\lambda^{m-1}S_{2,\lambda}(n,m),\quad (n\ge 1).
	\end{displaymath}
\end{theorem}
From \eqref{16}, we note that
\begin{align}
	t=l_{1,\lambda}\big(1-e_{\lambda}(-t)\big)&=\sum_{m=1}^{\infty}\frac{(1)_{m,1/\lambda}(-\lambda)^{m-1}}{m!}\big(1-e_{\lambda}(-t)\big)^{m} \nonumber \\
	&=\sum_{m=1}^{\infty}(1)_{m,1/\lambda}(-\lambda)^{m-1}\sum_{n=m}^{\infty}S_{2,\lambda}(n,m)(-1)^{n-m}\frac{t^{n}}{n!} \label{36} \\
	&=\sum_{n=1}^{\infty}\bigg(\sum_{m=1}^{n}(1)_{m,1/\lambda}\lambda^{m-1}(-1)^{n-1}S_{2,\lambda}(n,m)\bigg)\frac{t^{n}}{n!}. \nonumber
\end{align}
By comparing the coefficients on both sides of \eqref{36}, we obtain the following theorem.
\begin{theorem}
	For $n\in\mathbb{N}$, we have
	\begin{displaymath}
		(-1)^{n-1}\sum_{m=1}^{n}(1)_{m,1/\lambda}\lambda^{m-1}S_{2,\lambda}(n,m)=\delta_{n,1},
	\end{displaymath}
	where $\delta_{n,k}$ is Kronecker's symbol.
\end{theorem}
\emph{Remark.} Note that
\begin{displaymath}
	\lim_{\lambda\rightarrow 0}\beta_{n,\lambda}^{(1)}=(-1)^{n}B_{n},\quad \lim_{\lambda\rightarrow 0}\beta_{n,\lambda}^{(1)}(x)=(-1)^nB_{n}(x).
\end{displaymath}
From Theorem 8 and Theorem 9, we note that
\begin{displaymath}
	B_{0}=1,\quad B_{n}(1)-B_{n}=\left\{\begin{array}{cc}
		1, & \textrm{if $n=1$}, \\
		0 & \textrm{otherwise}.
	\end{array}\right.
\end{displaymath}
\begin{corollary}
	For $n\in\mathbb{N}$, we have
	\begin{displaymath}
		\sum_{m=1}^{n}(-1)^{n-m}(m-1)!S_{2}(n,m)=\left\{\begin{array}{cc}
		1, & \textrm{if $n=1$}, \\
		0 & \textrm{otherwise}.
	\end{array}\right.
	\end{displaymath}
\end{corollary}
It is well known that the Stirling numbers of the first kind are given by
\begin{equation}
	\frac{1}{k!}\big(\log(1+t)\big)^{k}=\sum_{n=k}^{\infty}S_{1}(n,k)\frac{t^{n}}{n!},\quad(\mathrm{see}\ [3,12,22]), \label{37}
\end{equation}
where $k$ is a nonnegative integer. \\
In view of \eqref{37}, we may consider the degenerate Stirling numbers of the first kind which are given by
\begin{equation}
	\frac{1}{k!}\big(\log_{\lambda}(1+t)\big)^{k}=\sum_{n=k}^{\infty}S_{1,\lambda}(n,k)\frac{t^{n}}{n!},\quad (k\ge 0). \label{38}
\end{equation}
From \eqref{38}, we note that
\begin{align}
	\sum_{k=0}^{\infty}\frac{(x)_{k,\lambda}}{k!}\big(\log_{\lambda}(1+t)\big)^{k}&=\sum_{k=0}^{\infty}(x)_{k,\lambda}\sum_{n=k}^{\infty}S_{1,\lambda}(n,k)\frac{t^{n}}{n!} \label{39}\\
	&=\sum_{n=0}^{\infty}\bigg(\sum_{k=0}^{n}S_{1,\lambda}(n,k)(x)_{k,\lambda}\bigg)\frac{t^{n}}{n!} \nonumber.
\end{align}
On the other hand,
\begin{align}
	\sum_{k=0}^{\infty}\frac{(x)_{k,\lambda}}{k!}\big(\log_{\lambda}(1+t)\big)^{k}&=e_{\lambda}^{x}\big(\log_{\lambda}(1+t)\big)=(1+t)^{x} \label{40} \\
	&=\sum_{n=0}^{\infty}(x)_{n}\frac{t^{n}}{n!}. \nonumber
\end{align}
By \eqref{39} and \eqref{40}, we see that the degenerate Stirling numbers of the first kind are also given by
\begin{equation}
	(x)_{n}=\sum_{k=0}^{n}S_{1,\lambda}(n,k)(x)_{k,\lambda},\quad (n\ge 0). \label{41}
\end{equation}
Note that
\begin{displaymath}
	\lim_{\lambda\rightarrow 0}S_{1,\lambda}(n,k)=S_{1}(n,k),\quad (n,k\ge 0).
\end{displaymath}
Now, we observe that
\begin{align}
	(x)_{n+1}&=(x)_{n}(x-n)=\sum_{k=0}^{n}S_{1,\lambda}(n,k)(x)_{k,\lambda}(x-k\lambda+k\lambda)-n\sum_{k=0}^{n}S_{1,\lambda}(n,k)(x)_{k,\lambda}\nonumber \\
	&=\sum_{k=0}^{n}S_{1,\lambda}(n,k)(x)_{k+1,\lambda}+\lambda\sum_{k=0}^{n}kS_{1,\lambda}(n,k)(x)_{k,\lambda}-n\sum_{k=0}^{n}S_{1,\lambda}(n,k)(x)_{k,\lambda} \label{42} \\
	&=\sum_{k=1}^{n+1}S_{1,\lambda}(n,k-1)(x)_{k,\lambda}+\sum_{k=0}^{n}(\lambda k-n)S_{1,\lambda}(n,k)(x)_{k,\lambda}. \nonumber
\end{align}
On the other hand,
\begin{equation}
	(x)_{n+1}=\sum_{k=0}^{n+1}S_{1,\lambda}(n+1,k)(x)_{k,\lambda}.\label{43}
\end{equation}
Therefore, by \eqref{42} and \eqref{43} and with the usual convention that $S_{1,\lambda}(n,k)=0$, for $k>n$ or $k<0$, we obtain the following theorem.
\begin{theorem}
	For $\lambda\in\mathbb{N}$, and $0\le k\le n+1$, we have
	\begin{displaymath}
		S_{1,\lambda}(n+1,k)=S_{1,\lambda}(n,k-1)+(\lambda k-n)S_{1,\lambda}(n,k).
	\end{displaymath}
\end{theorem}
From \eqref{18}, we note that
\begin{equation}
	\sum_{n=0}^{\infty}\beta_{n,\lambda}^{(k)}\frac{t^{n}}{n!}=\sum_{n=0}^{\infty}\frac{(-\lambda)^{n}(1)_{n+1,1/\lambda}}{n!(n+1)^{k}}\big(1-e_{\lambda}(-t)\big)^{n}.\label{44}
\end{equation}
Thus, by replacing $t$ by $-\log_{\lambda}(1-t)$, we get
\begin{equation}
	\sum_{m=0}^{\infty}\beta_{m,\lambda}^{(k)}\frac{(-1)^{m}}{m!}\big(\log_{\lambda}(1-t)\big)^{m}=\sum_{n=0}^{\infty}\frac{(-\lambda)^{n}(1)_{n+1,1/\lambda}}{n!(n+1)^{k}}t^{n}.\label{45}
\end{equation}
On the other hand,
\begin{align}
		\sum_{m=0}^{\infty}\beta_{m,\lambda}^{(k)}\frac{(-1)^{m}}{m!}\big(\log_{\lambda}(1-t)\big)^{m}\ &=\ \sum_{m=0}^{\infty}\beta_{m,\lambda}^{(k)}(-1)^{m}\sum_{n=m}^{\infty}S_{1,\lambda}(n,m)(-1)^{n}\frac{t^{n}}{n!}\label{46}\\
		&=\ \sum_{n=0}^{\infty}\bigg(\sum_{m=0}^{n}(-1)^{n-m}\beta_{m,\lambda}^{(k)}S_{1,\lambda}(n,m)\bigg)\frac{t^{n}}{n!}. \nonumber
\end{align}
Therefore, by \eqref{45} and \eqref{46}, we obtain the following theorem.
\begin{theorem}
	For $k\in\mathbb{Z}$ and $n\ge 0$, we have
	\begin{displaymath}
		\frac{1}{(n+1)^{k}}=\frac{1}{\lambda^{n}(1)_{n+1,1/\lambda}}\sum_{m=0}^{n}(-1)^{m}\beta_{m,\lambda}^{(k)}S_{1,\lambda}(n,m).
	\end{displaymath}
\end{theorem}

\end{document}